\documentclass{article}
\usepackage[utf8]{inputenc}
\usepackage[tbtags]{amsmath}
\usepackage{amsfonts,amssymb,mathrsfs,amscd}

\def \Z {{\mathbf {Z}}}

\def \N {{\mathbb {N}}}

\def\uu{\bigsqcup}

\def\s{ S}

\begin{document}

\title{\LARGE    \bf Spectra and joint dynamics of Poisson suspensions \\
for rank-one automorphisms}
\author{\LARGE   V.\,V.~Ryzhikov}
%\address{Moscow State University}
%\email{vryzh@mail.ru}

\date{}
%\udk{517.987}

\maketitle

%\begin{fulltext}

\Large

\begin{abstract}{For every natural $n>1$, there is unitary operator  $T$ of dynamical origin such that its tensor power $T^{\otimes n}$ has singular spectrum, and $T^{\otimes (n+1)}$ has absolutely continuous one. For a set $D$ of positive measure there are mixing zero entropy automorphisms $S,T$ such that $S^nD\cap T^nD=\varnothing$ for all $n>0$. The following answers, in particular, Frantzikinakis-Host's question. If $ p(n+1)- p(n), \ \ q(n+1)- q(n)\ \to\ +\infty$, then there is a divergent sequence $ \sum_{n=1}^{N} \mu(S^{ p(n)}C\cap T^{ q(n)}C)/N$ for some set $C$ and automorphisms $S,T$ of simple singular spectrum.

Key words: \it Sidon rank one constructions, tensor products, spectrum, Poisson suspensions, joint dynamics, divergence of ergodic averages.\rm}
\end{abstract}
%\begin{keywords}

%\end{keywords}
\markright{Spectra and joint dynamics of Poisson suspensions}

%%%%%%%%%%%%%%%%%%%%%%%%%%%%%%%%%%%%%%%%%%%%%%%%%%%%%%%%%%%%%%%%%%%%%%%%%%%%%%%%%%%%%%
\section{Introduction}
By the spectral theorem, a unitary operator $S$ on a separable Hilbert space is isomorphic to the multiplication operator $V$ acting
in the space
$L_2(\mathbb T\times \mathbb N, \sigma)$ by the formula $Vf(t,n)=tf(t,n).$ Here $\mathbb T$ is the unit circle in the complex plane, and $\sigma$ is the normalized Borel measure on $ \mathbb T\times \mathbb N$.
If the measure $\sigma$ is singular, then the spectrum of the operator $S$ is said to be singular. If the projection of the measure $\sigma$ onto $\mathbb T$ is absolutely continuous with respect to the Lebesgue measure on the circle, then its spectrum is called absolutely continuous.

\vspace{2mm}
\bf Theorem 1. \it For every natural $n>1$ there exists a unitary operator $T$ such that its tensor power $T^{\otimes n}$ has  singular spectrum, and  spectrum of power $T^{\otimes (n+1)}$ is absolutely continuous. Thus, the convolution power $\sigma^{\ast n}$ of the spectral measure $\sigma$ of such an operator $T$ is singular, and $\sigma^{\ast (n+1)}$ is absolutely continuous.\rm

\vspace{2mm}
 The proposed operator $T$ has a dynamical origin. It is induced by a Sidon automorphism of a space with sigma-finite measure, see \cite{24msb}. In the paper, automorphisms and the operators induced by them are denoted identically.

In \cite{F}, a question was asked about the convergence of ergodic averages for a sequence $S^{n^2}f\cdot T^{n^3}g$, where  zero entropy $S,T$  are ergodic automorphisms of a probability space , and $f,g$ are
bounded measurable functions. In this connection, a general problem arose
of finding  divergent averages
$$ \frac 1 N \sum_{n=1}^N \int S^{ p(n)}f \cdot T^{ q(n)}g \, d\mu. $$  for some zero entropy automorphisms  $S,T$  and sequences ${ p (n)},{ q (n)}\to\infty$. 

\vspace{3mm}
\bf Theorem 2. \it
Let $ p(n+1)- p(n)$, $ q(n+1)- q(n)\ \to\ +\infty.$
There exist automorphisms $S,T$ of simple singular spectrum such that the sequence $$ \frac 1 N \sum_{n=1}^{N} \mu(S^{ p(n)}C\cap T^{ q(n)}C)$$ diverges for some measurable set $C$. \rm

\vspace{2mm}
The divergence of ergodic averages in the case $p(n)=q(n)=n$ for zero entropy automorphisms was discovered by different methods in \cite{A} and \cite{Ye}. In \cite{R}, Poisson suspensions were used for a similar result. Theorem 2 strengthens some results of \cite{Ye24}, \cite{K}. It is proved by use of Poisson suspensions over\footnote{"Poisson suspensions" = in \cite{KSF} are literally written as "Poisson  superconstructions", so we write "over".} specially chosen rank-one automorphisms.
The following statement  answers  Frantzikinakis' question.

\vspace{2mm}
\bf Theorem 3. \it
For some mixing zero entropy automorphisms $S,T$   and a set $D$ of positive measure, $S^nD\cap T^nD=\varnothing$ holds for all $n>0$. \rm

\vspace{2mm}
As suitable $S,T$, we choose tensor powers of Poisson suspensions  over Sidon automorphisms. Infinite automorphisms (automorphisms of a space with sigma-finite measure) turned out to be a convenient tool in solving problems on the joint dynamics of automorphisms of a probability space.
\section{ Rank one  Sidon constructions}
The proofs of the proposed theorems use well-known rank one constructions  in ergodic theory. The first applications were given by R.V. Chacon \cite{Ch} and D. Ornstein \cite{Or}.  D. Rudolph discovered the property of minimal self-joinings for a mixing automorphism of rank one, the base of an impressive machine of counterexamples \cite{Ru}. 

Among modern examples of the application of rank one constructions,  mention an application to  V.I. Oseledets' problem. There is a measure on a line with complete support whose tensor square has absolutely continuous projections onto a fixed line in the plane along one dense set of directions and singular projections along another dense set \cite{R22}. Such measures turned out to be spectral measures of some constructions of rank one. Constructions of rank one were used to solve the problem of \\ I.V. Podvigin on typicality in the Alpern-Tikhonov space of slow decays of correlations for mixing automorphisms \cite{mz24}. Staircase rank one constructions  underlie the effective solution of A.N. Kolmogorov's problem on the group property of  spectrum and V.A. Rokhlin's problem on  homogeneous spectrum of multiplicity 2 in the class of mixing automorphisms \cite{25}. Applications of Poisson suspensions over rank one automorphisms will be discussed later.

\bf  Rank one onstructions. \rm Fix a natural number $h_1$, a sequence $r_j$ (the number of columns into which the tower of stage $j$ is cut) and a sequence of integer vectors (spacer parameters)
$$ \bar s_j=(s_j(1), s_j(2),\dots, s_j(r_j-1),s_j(r_j)).$$
Let at step $j$ a system of non-intersecting half-intervals
$$E_j, TE_j, T^2E_j,\dots, T^{h_j-1}E_j,$$ be defined, and on the half-intervals $E_j, TE_j, \dots, T^{h_j-2}E_j$
the transformation $T$ is a parallel translation. Such a set of half-intervals is called the tower of stage $j$, their union is denoted by $X_j$ and is also called a tower.

%Let us represent $E_j$ as a disjoint union of $r_j$ half-interval$E_j^1,E_j^2E_j^3,$ $\dots E_j^{r_j}$ of the same length. For each $i=1,2,\dots, r_j$, we define a column $X_{i,j}$ as the union of intervals
$E_j^i, TE_j^i ,T^2 E_j^i,\dots, T^{h_j-1}E_j^i.$
To each column $X_{i,j}$, we add $s_j(i)$ disjoint half-intervals of the same measure as $E_j^i$, obtaining a set
$$E_j^i, TE_j^i, T^2 E_j^i,\dots, T^{h_j-1}E_j^i,
T^{h_j}E_j^i, T^{h_j+1}E_j^i, \dots, T^{h_j+s_j(i)-1}E_j^i$$
(all these sets are disjoint).
Denoting $E_{j+1}= E^1_j$, we set
$$T^{h_j+s_j(i)}E_j^i = E_j^{i+1}, \ i<r_j.$$
The added columns line up to form a tower at stage $j+1$:
$$E_{j+1}, TE_{j+1}, T^2 E_{j+1},\dots, T^{h_{j+1}-1}E_{j+1},$$
where $ h_{j+1} =h_jr_j +s_j(1)+s_j(2)+\dots +s_j(r_j).$
The partial definition of the transformation $T$ at stage $j$ is preserved at all subsequent stages. As a result, on the space $X=\cup_j X_j$, an invertible transformation $T:X\to X$ is defined that preserves the standard Lebesgue measure on $X$. Rank one automorphism  is ergodic, it has  simple  spectrum.
It is also known that linear combinations of interval (floor) indicators appearing in the description of the constructions of $T$ are cyclic vectors for the operator $T$ with  continuous spectrum. In the case of infinite measure  space $X$, every ergodic automorphism of this space has  continuous spectrum.

\bf Sidon automorphisms. \rm
Let a rank one construction $T$  have the following property:
\it the intersection
$T^mX_j$ for $h_{j}<m\leq h_{j+1}$ can be contained
only in one of the columns $X_{i,j}$ of the tower $X_j$. \rm Such a construction  is called \it Sidon. \rm
\vspace{2mm}

\bf Sidon automorphisms of class $C(\nu)$. \rm
We say that a Sidon construction $T$ of rank one
belongs \it to class $C(\nu)$, $\nu\geq 0$ \rm if for some sequence
$j(k)\to\infty$ for any $\delta>0$ the following holds:
$$\sum_k r_{j(k)}^{-\delta}<\infty;$$
$$r_j=r_{j(k)}, \ j(k)\leq j<j(k+1); $$
$$j(k+1) -j(k)=[r_{j(k)}^\nu].$$
For convenience, we will assume that the condition
$$ h_j\ll s_j(1)\ll s_j(2)\ll \dots\ll s_j(r_j-1)\ll s_j(r_j), \ \ r_j\to\infty,$$
is satisfied, where expressions of the form $a_j\ll b_j$ mean that $b_j>\psi(j)a_j$ for some fixed sequence $\psi(j)\to +\infty$.

\vspace{2mm}
\bf Theorem 4. \it Let $T\in C(\nu)$ be a Sidon automorphism, $\nu\geq 0$.
For $\nu\geq 2n-2$, its power  $T^{\otimes n}$ has  singular spectrum, and for $\nu< 2n-2$ its spectrum is absolutely continuous.
\rm

\vspace{2mm}

Conservativity of automorphism $S$ means the following: if $\mu(A)>0$, then there exists $n>0$ such that $\mu(A\cap S^nA)>0$. Dissipativity of automorphism $T$ is defined as follows: for some measurable set $W$, all sets $T^nW$, $n\in \Z$, are disjoint, and their union gives the entire space.
It should be said that 

\it for $T\in C(\nu)$ for $\nu\geq n-1$ the power $T^{\otimes n}$ is conservative, and for $\nu < n-1$ it is dissipative. \rm

This follows from a more general fact: for a Sidon construction $T$, the fact that
$$\sum_{i=1}^\infty\left(\frac 1 {r_j}\right)^{d-1} =\infty$$
is equivalent to the fact that  $T^{\otimes d}$ is conservative. This statement is a consequence of Theorem 3.2 of \cite{LS}.

\section{ Singular  spectrum of  $T^{\otimes d}$}
We will show that for the Sidon construction $T\in C(\nu)$ the spectrum
of power $T^{\otimes d}$ is singular if the condition
$$\sum_{k=1}^{\infty }\left(\frac 1 {r_{j(k)}}\right)^{2d-2-\nu}=\infty,$$ \rm
 the case when $\nu>2d-2$.

Let $F$ be the indicator of the Cartesian $d$-cube of the tower $A=X_1$, and $S$ denote the restriction of power $T^{\otimes d}$ to the cyclic space of the vector $F$.
Since $\chi_{A}$ is a cyclic vector for the operator $T$,
the spectral types of the operators $S$ and $T^{\otimes d}$ are the same. We will prove the singularity of the spectrum of the operator $S$.

 We set $q(1,j)=0$ and $$ q(i,j)=(i-1)h_j +s_j(1)+s_j(2)+\dots +s_j(i-1),\, 1<i\leq r_j,$$
$$Q_j(S)= \sum_{1\leq i\neq i'\leq r_j} S^{q(i,j)-q(i',j)}.$$
Denote $ a_k= (r_{j(k)}-1)r_{j(k)}^{1-d}, \ \ N_k=j(k+1)-j(k)= \left[r_{j(k)}^\nu\right].$ Our goal is to show for the sequence of polynomials
$$ P_k = (a_kN_k)^{-1}\sum_{j=j(k)+1}^{j(k+1)} Q_j$$
the convergence of $$P_k(S)F\to F, \ \ P_k(S^p)F\to 0, \ p>1.$$

By Lemma 2.1 \cite{24msb} these convergences guarantee the singularity of the spectrum
of the operator $S$. Let us briefly recall the proof of this lemma.
Let a bounded operator $J\neq 0$ intertwine $S$ and its power $S^p$. Note that
in our case $P_k(S)F\to F$ is equivalent to $P_k(S^\ast)F\to F$. Then for arbitrary linear combinations $V, W $ of shifts $S^i F$ of the cyclic vector we have
$$(JV,W) = \lim_k(JP_k(S^\ast)V,W)=\lim_k(JV,P_k(S^p)W)= 0.$$
But this means that $J=0$, hence the operators $S$ and $S^p$, $p>1$, are spectrally
disjoint. In this case, the spectrum of the operator $S$ is singular.

Consider the functions 
$$ G_k =\sum_{j=j(k)}^{j(k+1)-1}Q_j(S)F. $$ 
Let's represent $G_k$ in the form 
$$G_k= a_kN_kF +(G_k-a_kN_k)F +G_k(1-F), \ a_k= (r_j-1) r_j^{1-d}.$$ 
Let's make sure that for $\nu>2d-2$ 
$$ \frac {\left\|(G_k-a_kN_k)F  + G_k(1-F)\right\|} {\|a_kN_kF\|} \to 0,  \ k\to\infty. \eqno (3.1)$$

\bf Norm estimate of  $\bf (G_k-a_kN_k)F$. \rm
Above, we have already defined the value of $a_k$ without explanation. Now let's explain its origin: $$a_k=\int_{X_1^{\times d}} D_{j(k)} \, d\mu^d, \ \ D_{j}= \sum_{1\leq i\neq i'\leq r_j} F \cdot S^{q(i,j)-q(i',j)}F.$$ On the union $\uu_{i=1}^{r_j } E_{i,j}^{\times d}$ the function $D_j$ is equal to $r_j-1$, and outside this set it is equal to $0$, so $$ \|D_j\|^2= (r_{j(k)}-1)^2 r_{j(k)}^{1-d}< r_{j(k)}^{3-d}, \ j(k)\leq j<j(k+1).$$

Important property: \bf functions $D_j$ as random variables on the probability space $(X_1^{\times d},\mu^d)$ are independent. \rm 

This follows directly from the independence in $(X_1,\mu)$ of the sets $X_{i,j}\cap X_1$ and $X_{i',j'}\cap X_1$ for $j\neq j'$. If we subtract the constant $a_k=(r_j-1)r_{j(k)}^{1-d}$ (the value of their integral) from the functions $D_j$,
then the resulting differences $\Delta_j$ will be pairwise orthogonal in $L_2(X_1^{\times d}, \mu^d)$ and have the same norm in the case $j(k)\leq j <j(k+1)$.
The square of the norm of their sum of such $\Delta_j$ is $\|(G_k-a_kN_k)F\|^2$. So as $(a_kN_k)^{-2}= r_{j(k)}^{-4 +2d -2\nu}$, we have $$\|a_kN_kF\|^{-2}\|(G_k-a_kN_k)F\|^2\leq r_{j(k)}^{-4 +2d -2\nu +3-d+\nu}=
r_{j(k)}^{d-1-\nu}\to 0, \ k\to\infty.$$

\bf Norm estimate of $\bf (1-F)G_k$. \rm By definition
$$(1-F)G_k= \sum_{j=j(k)}^{j(k+1)-1}(1-F)Q_j(S)F.$$
The support of the function $(1-F)Q_j(S)F$ lies in $X_{j+1}\setminus X_j$, therefore the supports
of the functions $(1-F)Q_j(S)F$ do not intersect.
Let us show that the function $(1-F)Q_j(S)F$ is an indicator of its support, the measure of which is less than $r_j^2-r_j$.

\vspace{2mm}
\bf Lemma. \it For $(i,i')\neq (m,m')$   the sets
$T^{q(i,j)-q(i',j)}X_1$ and $T^{q(m,j)-q(m',j)}X_1$ utside $X_1$ are disjoint. \rm

\vspace{2mm}
This lemma was implicitly proved in \cite{24msb} in the proof of item (ii) of Theorem 4.1.
Let us clarify the main points. For our constructions, the mapping $T^{q(i,j)-q(i',j)}$ sends the $i'$-th column of $X_{i',j}$ to the $i$-th column of $X_{i,j}$, and the images of the remaining columns for $i\neq i'$ end up in $X_{j+1}\setminus X_j$. If
$$\mu(T^{q(i,j)-q(i',j)}X_1\cap T^{q(m,j)-q(m',j)}X_1)>0,$$
this means that
$q(i,j)-q(i',j)-q(m,j)+q(m',j)$ coincides (!) with $q(n,j)-q(n',j).$ Since $h_j\ll q(1,j)\ll q(2,j) \dots \ll q(r_j-1,j),$
we get that the specified coincidence (!) is fulfilled only if $n,n'$ are among $i,i'm,m'$.
But this entails the coincidence of two numbers from the set $i,i'm,m'$. Let, for example, $i=m$,  and $i'\neq m'$.
Then the intersection
$T^{q(i,j)-q(i',j)}X_1 \cap T^{q(i,j)-q(m',j)}X_1$
contains $ X_{i,j}\subset X_1$, that was required to be shown.

It follows from the lemma that outside $X_1^{\times d}$
the sets $S^{q(i,j)-q(i',j)}X_1^{\times d}$ do not intersect (recall that
$S= T^{\otimes d}$). Then
$$\|G_k(1-F)\|^2=\mu^n (supp\, G_k(1-F)) < (r^2_{j(k)}-r_{j(k)})N_k.$$
We get
$$\|a_kN_kF\|^{-2}\|G_k(1-F)\|^2
< r_{j(k)}^{d-2 -\nu}\to 0, \ k\to\infty, $$
this yields $(3.1)$. We have proved that $P_k(S)F\to F$.

The convergence of $P_k(S^p)F\to 0$ is a consequence of the fact that all functions $S^{pq(i,j)-pq(i',j)}F$ have disjoint supports (the corresponding argument is also available in \cite{24msb}). This leads to an estimate similar to the one above:
$$\|P_k(S^p)F\|^2 = (a_kN_k)^{-2} (r_{j(k)}^2 -r_{j(k)})N_k \ < const \cdot r_{j(k)}^{ 2d -2 -\nu} \to 0, \ k\to\infty.$$

\bf The case $\nu=d-2$. \rm We will show that for $\nu=d-2$ the tensor power $S$ has  singular spectrum. Now, in contrast to (3.1), we have
$$ \frac {\|(G_k-a_kN_k)F\|} {\|a_kN_kF\|} \to 0, \ \ k\to\infty, \ \
\frac {\| G_k(1-F)\|} {\|a_kN_kF\|} \ <\ Const. $$
But the supports of the functions $G_k(1-F) $ do not intersect, so repeated averaging  gives the desired result:
for example,
$$R_N(S)F=\frac 1 N \sum_{k=N+1}^{2N} P_k(S)F\ \to \ F , \ N\to\infty.$$
Moreover, for $p>1$, this is obvious, $R_N(S^p)F\to 0$, $N\to\infty$.
The spectrum of the operator $S$ is singular.

Thus, for $\nu \geq 2d-2$  spectrum of 
$T^{\otimes d}$ is singular.

\section{\bf Absolutely continuous  spectrum of $T^{\otimes d}$}
\bf Lemma. \it If for a Sidon construction $T$ the series
$$\sum_{i=1}^{\infty }\left(\frac 1 {r_j}\right)^{2d-2},$$ converges
then for $f$, the indicator of the tower $X_1$, we have
$\sum_{n}^\infty \left(T^nf\,,\,f\right)^{2d}\, <\, \infty,$
so $T^{\otimes d}$ has absolutely continuous spectrum.\rm

\vspace{2mm}
Proof. By induction on $m$ we establish the equality $$\sum_{|n|<h_{m}}\,\, \left(T^nf\,,\,f\right)^{2d} = \prod_{j=1}^{m-1} \left( 1+\frac{r_j^2-r_j} {r_j^{2d}} \right). \eqno (4.1)$$ For $h_{j} <|n|<h_{j+1}$ $\left(T^nf\,,\,f\right)=1/r_j$ only in those cases when the column $X_{i,j}$ under the action of $T^n$ moves into the column $X_{i',j}$, $i\neq i'$. Here it is implied that the lower floor of column $X_{i,j}$ will be at the place of the lower floor of column $X_{i',j}$.
The number of such cases is equal to ${r_j^2-r_m}$.
Let $c_m= \mu(T^m X_1\cap X_1)$, $|m|<h_j$, and $(T^nf\,,\,f)=1/r_j$, $h_{j}
<|n|<h_{j+1}$.
Then $$(T^{m+n} f\,,\,f)=c_m/r_j,$$
and all values of $c_p\neq 0$ for
$h_{j}<|p|<h_{j+1}$ take place only for $p=m+n$ for $m,n$ satisfying the above conditions. From the above we obtain that
$$\sum_{|p|<h_{j}}\, c_p^{2d} =\left( 1+\frac{r_j^2-r_j} {r_j^{2d}} \right)
\sum_{|m|<h_{j-1}}\, c_m^{2d},$$
this yields (4.1). The infinite product 
$$\prod_{j=1}^{\infty} ( 1+ ({r_j^2-r_j}) {r_j^{-2d}})$$ 
converges if
$$\sum_{j}({r_j^2-r_j}) {r_j^{-2d}} \ <\ \infty. $$
This series converges for $T\in C(\nu)$ for $\nu < 2d-2$, what we need.  Theorems 1 and 4 are proved.

\vspace{2mm}
\bf Remark. \rm For the Sidon automorphism $T\in C(1)$, we have
$$\sum_{n}^\infty \mu(T^n X_1\cap X_1)^{4}\, <\, \infty.$$
We will use this fact in the proof of Theorem 3.

Thus, for $T\in C(2)$, we obtain that $T$, $T^{\otimes 2}$ have  singular spectrum, the conservative power $T^{\otimes 3}$ has  absolutely continuous spectrum,
and $T^{\otimes 4}$ and the remaining powers are dissipative (their spectrum is countably multiple Lebesgue).
For $T\in C(5)$ the spectrum of powers $T$, $T^{\otimes 2}$ and $T^{\otimes 3}$ is singular, $T^{\otimes 4}$, $T^{\otimes 5}$, $ T^{\otimes 6}$ are conservative, their spectrum is absolutely continuous, the remaining powers are dissipative.

It is likely that $T^{\otimes 3}$ also has a countably multiple Lebesgue spectrum,
but proving this in the general case requires effort. In special cases,
when $T\in C(5)$, by controlling the parameters of the constructions, it is easy to ensure that the spectrum of $T^{\otimes 3}$ is Lebesgue.

Besides the Banach problem about the existence of an automorphism with
 simple Lebesgue spectrum (automorphisms of rank one can still claim to have a Lebesgue spectrum), there is perhaps an even more difficult problem in the spectral theory of transformations:  \it is there a measure space  automorphism  with  absolutely continuous but non-Lebesgue spectrum? \rm

\section{Poisson suspensions}
Bernoulli automorphisms with different entropy form a continuous family of pairwise non-isomorphic automorphisms with the same Lebesgue spectrum, that  was established by A. N. Kolmogorov \cite{Ko}.
 It became possible to find  families  of zero entropy  non-isomorphic automorphisms having the same spectrum thanks to Poisson suspensions. The latters are  known both in the representation theory of large groups \cite{VGG}-\cite{Ne} and in the ergodic theory \cite{KSF}-\cite{PR}. In the paper \cite{PR}, F. Parreau and E. Roy proposed a large variety of families of non-isomorphic Poisson suspensions with the same singular spectrum.  Let us give other examples  of Poisson suspension applications.

M.I. Gordin asked the author the following question: can automorphisms with singular spectrum have an ergodic homoclinic group?
 In solving this problem, suitable examples turned out to be mixing Gaussian and Poisson suspensions over some automorphisms of rank one (see references in \cite{R23}). To prove Theorem 3 of this paper, we will make  use of homoclinic elements of the Poisson suspension over  rank one automorphisms.

Bergelson's question about rigid and mixing sequences of automorphism powers was answered by using Poisson suspensions over automorphisms of rank one. In \cite{R21} it is shown that for every set $M$ of zero density there exists a rigid transformation ($T^{n_i}\to I, \ n_i\to\infty$) mixing along the set $M$ ($T^{m}\to_w 0$ in the space $Const^\perp$ for $m\in M$ and $ m\to\infty$).

The Poisson suspensions over rank-one automorphisms gives all possible sets of spectral multiplicities of the form $M\cup \{\infty\}$ ($M\subset \N$) \cite{R23},\cite{24msb} (weakly mixing and mixing).

Examples of zero entropy automorphisms with  completely positive $P$-entropy (Kirillov-Kushnirenko invariant) were found again among Poisson suspensions over rank one automorphisms \cite{RT}.

\bf Poisson measure. \rm Consider the configuration space $X_\circ$ consisting of all infinite countable sets $x_\circ$ such that each interval of the spaces $X$ contains only a finite number of elements of the set $x_\circ$.

The space $X_\circ$ is equipped with the Poisson measure. Recall its definition.
Subsets $A\subset X$ of finite $\mu$-measure in the configuration space $X_\circ$ are assigned cylindrical sets
$C(A,k)$, $k=0,1,2,\dots$, by the formula
$$C(A,k)=\{x_\circ\in X_\circ \ : \ |x_\circ\cap A|=k\}.$$

All possible finite intersections of the form $\cap_{i=1}^N C(A_i,k_i)$
form a semiring.
On this semiring we define the measure $\mu_\circ$ as follows.
Under the condition that measurable sets $A_1, A_2,\dots, A_N$ are disjoint
and have finite measure, we set
$$\mu_\circ(\bigcap_{i=1}^N C(A_i,k_i))=\prod_{i=1}^N \frac {\mu(A_i)^{k_i}}{k_i!} e ^{-\mu(A_i)}.$$
Continuing the measure from the semiring of cylinders leads to the Poisson space
$(X_\circ,\mu_\circ)$, which is isomorphic to the standard Lebesgue probability space. An automorphism $T$ of $(X,\mu)$ naturally induces an automorphism
$\mathring T$ of $(X_\circ,\mu_\circ)$, which is called the Poisson suspension of $T$.

As an operator, $\mathring T$ is isomorphic to the operator
$$\exp(T)= {\bf 1} \oplus T \oplus T^{\odot 2}\oplus T^{\odot 3}
\oplus\dots .$$ Note that the Gaussian automorphism induced by the orthogonal operator $T$ is also isomorphic to $\exp(T)$, but its metric properties may differ greatly from the metric properties of the suspension $\mathring T$ (see \cite{PR}). Thanks to the above results on Sidon automorphisms, we get new spectral types of Poisson  suspensions of the form
$$\sigma+\sigma^{\ast 2}+\dots+ \sigma^{\ast n}+ \lambda, \ n>1,$$
where the convolution powers $\sigma^{\ast n}$ are singular, and $\lambda$ is Lebesgue. All terms in the above sum are mutually singular.
 
\section{ Proof of Theorem 2}

Given sequences ${ p(n)},{ q(n)}$, where
$ p(n+1)- p(n),$ $q(n+1)- q(n)\ \to\ +\infty. $
It is required to find a set $C$ and Poisson suspensions $S,T$ with simple singular spectrum such that the sequence $ \sum_{n=1}^{N} \mu(S^{ p(n)}C\cap T^{ q(n)}C)/N$ diverges.

\bf Definition of automorphism $\bf S$. \rm Consider an infinite rank one automorphism $\s$ , defined by parameters
$r_j=2$, $ s_j(1)=0$, $ s_j(2)\to \infty$. In this case, the parameters $s_j(2)$ are chosen so that $s_j(2)>\max\{p(jh_j),q(jh_j)\}$. If sets $A,B$ consist of floors of the tower of stage $j_0$,
then for $j>j_0$ we have
$$\mu(\s^{h_j}A\cap B)=\mu(A\cap B)/2.$$
Therefore $ \s^{h_j}\to_w I/2$, and this implies the singularity of the spectrum of the operator
$$\exp(\s)= {\bf 1} \oplus \s\oplus\s^{\odot 2}\oplus\s^{\odot 3}\oplus\dots \ $$
Therefore, the Gaussian and Poisson suspensions over $\s$ (isomorphic to the operator $\exp(\s)$) have singular spectra. Their entropy is zero.
Later, we modify the construction of $\s$, ensuring simple 
 spectrum of  $\exp(\s)$.

\bf Definition of automorphism $\bf T$. \rm Automorphism $\s$ cyclically permutes the floors in tower $X_j$ (except the top floor).
Automorphism $T$ will permute these floors in a different order. We define it as
$P S P^{-1}$, where $P$ also permutes floors.

Let us identify the set of floors of the tower of stage $j+1$ with the ordered
set of numbers $I_j=(0,1,2,\dots, h_{j+1}-1)$. Automorphism $\s$
obviously corresponds to an (almost) cyclic permutation $\sigma$ on $I_j$. It is given by the formula $\sigma(n)=n+1$, but let the value $\sigma(h_{j+1}-1)$ be undefined. (In the spirit of  Katok-Stepin's approximations, we can formally set $\sigma(h_{j+1}-1)=0$.)

For $n>n_j$ intervals $(p(n),p(n)+1,p(n)+2, \dots, p(n)+2h_j-1)$
do not intersect, intervals
$(q(n)+2h_j,q(n)+2h_j+1,q(n)+2h_j+2, \dots,q(n) +4h_j-1)$
do not intersect either.

 \bf For odd $\bf j$ \rm ensure $2h_j<p(n),q(n)< h_{j+1} - 4h_j$ and choose a permutation $\pi_j$ on $\{2h_j, 2h_j-1, \dots, h_{j+1}-1\}$ such that $$\pi_j(q(n)+i)= p(n)+2h_j +i, \ i=0,1, 2,\dots, 2h_j.$$
\\ With this permutation $\pi_j$ we associate an automorphism $P_j$ that sends the $k$-floor to the $\pi_j(k)$-floor by parallel translation.
Outside $X_{j+1}\setminus X_j$, let the permutation of floors $P_j$
act identically. Recall that $$X_{j}=\uu_{i=0}^{2h_{j}-1} \s^iE_{j+1}, \ \
X_{j+1}=\uu_{i=0}^{h_{j+1}-1} \s^iE_{j+1}.$$

The set $A=X_1\subset X_j$ consists of the floors of stage $j+1$
(we identify $A$ with the set of numbers of these floors, which form a subset of the set $\{0,1,2, \dots, 2h_j\}$.)
Set $T=P_j\s P_j^{-1}$, then
$$T^{q(n)}A= \s^{p(n)+2h_j}A,$$
$$T^{q(n)}A\cap \s^{p(n)}A=\varnothing, \eqno (6.1) $$

\bf For even $\bf j$ \rm define a permutation $\pi_j$:
$$\pi_j(q(n),q(n)+1, \dots, q(n)+2h_j-1)=
(p(n),p(n)+1, \dots,p(n) +2h_j-1),$$
$2h_j<p(n),q(n)< h_{j+1} - 2h_j$.
The corresponding automorphism $P_j$ again maps the $k$-th floor to
to the $\pi_j(k)$-th floor, and outside $X_{j+1}\setminus X_j$ we set $P_j=Id$.

For $T=P_j\s P_j^{-1}$ we have
$$T^{q(n)}A=\s^{p(n)}A.\eqno (6.2)$$

Let $$P=\prod_j P_j, \ \ T=P\s P^{-1}.$$ It is important that for such $T$ equalities (6.1), (6.2) are preserved.

So, we can find $N_j\to\infty$ (with corresponding $h_j\to\infty$) such that
$$ \frac 1 {N_{2k}} \sum_{n=1}^{N_{2k}} \mu(\s^{ p(n)}A\cap T^{q(n)}A)\to \mu(A), $$
and
$$ \frac 1 {N_{2k+1}} \sum_{n=1}^{N_{2k+1}} \mu(\s^{ p(n)}A\cap T^{ q(n)}A)\to 0.$$
We choose $N_j=\max \{n\,:\, p(n),q(n)<h_{j+1}-4h_j$.
Note that $N_{j-1}/N_{j}\to 0$, since $N_{j-1}<h_j$, $N_{j}>jh_j$ (recall that we chose $h_{j+1}>\max\{p(jh_j),q(jh_j)\}$ above).

\bf Modification. \rm For some rare set of stages $j$, we can modify the parameters of the above construction $\s$ so that the weak closure of its powers contains all possible polynomials. To do this, let $r_j=2j$, $s_j(i)=0$ for $i\leq j$, and $s_j(i)=1$ for $i\leq 2j$ (Katok's spacers). This gives a simple spectrum for the operator $\exp(\s)$. For these stages, let $P_j=Id$.

\bf Completion of the proof. \rm Consider in the Poisson space the cylinder $C=C(A,0)$, $0<c=\mu_\circ(C)<1$, $\mu_\circ$ is the Poisson measure. If the sets $A,A'$ do not intersect, then the corresponding cylinders $C,C'$ are independent
with respect to the measure $\mu_\circ$.
Since $\mathring S^n C(A,0)=C(\s^nA,0)$, $\mathring T^n C(A,0)=C(T^nA,0)$, we get $$ \frac 1 {N_{2k+1}} \sum_{n=1}^{N_{2k+1}} \mu_\circ(\mathring S^{p(n)}C\cap \mathring T ^{q(n)}A)\to c^2,$$ $$ \frac 1 {N_{2k}} \sum_{n=1}^{N_{2k}} \mu_\circ(\mathring S^{p(n)}C\cap \mathring T^{q(n)}A)\to c.$$

\section{Proof of Theorem 3} Now we need to present Poisson suspensions $\mathring S,\mathring T$ of zero entropy  such that for a set $D$ of positive measure $\mathring S^nD\cap \mathring T^nD=\varnothing$ for all $n>0$.
The Gordin homoclinic group $H(\mathring T)$ is defined as
$$H(\mathring T)=\{R\in Aut(\mu_\circ) \,:\, \mathring T^{-n}\mathring R \mathring T^n\, \to_s\, I, \ n\to\infty\}.$$
The mixing Poisson suspension over an infinite automorphism $T$ (it must satisfy the condition $T^n\to_w 0$) has the following homoclinic elements. Let $ \mu ( x: R x\neq x)<\infty$. Then $T^{-n}R T^n\, \to_s\, { I}$, so for the corresponding suspensions $\mathring R,\mathring T$ we have
$$\mathring T^{-n}\mathring R \mathring T^n\, \to_s\, {\ I}.$$ Now let $T$ be a rank-one transformation and $R\in [T]$ an involutive permutation of some disjoint floors $E$ and $R E$ in tower $X_1$. (Note that we can consider tower $X_j$ as $X_1$, and change the numbering of towers.) Outside $U:=E\cup R E$ let $R$ coincide with $Id$. In the Poisson space, consider the set $$D=C(E,0)\cap C(RE,1).$$ Recall that $C(E,0)$ consists of all countable sets $x_\circ=\{x_i\}$ such that $E\cap x_\circ$ has cardinality 0. (By definition, $x_\circ$ has no limit points in the intervals that make up the space $X$.)
The cylinder $C(R E,1)$ is formed by all $x_\circ$ such that $\mathring RE\cap x_\circ$ has cardinality 1. Therefore, $$\mathring R D\cap D=\varnothing.$$
We set $\mathring S=\mathring R^{-1}\mathring T\mathring R.$ Since $\mu_\circ (\mathring T^{-n}\mathring R\mathring T^n D\Delta D)\to 0,$ we get
$$\mu_\circ (\mathring S^nD\cap \mathring T^n D)= \mu_\circ (\mathring RD\cap \mathring T^{-n}RT^n D)\ \to 0.$$
If $$\sum_n \mu_\circ (\mathring RD\cap \mathring T^{-n}\mathring R\mathring T^n D)<\infty,$$ we achieve the goal as follows. For sufficiently large $p$ we take the power of $\mathring T^p$ instead of $\mathring T$. Since
$$\sum_n \mu_\circ (\mathring RD\cap \mathring T^{-pn}\mathring R\mathring T^{pn} D)\to 0,
\ p\to\infty,$$
there is a subset $D'\subset D$ of positive measure such that
$$\mathring S^{pn}D'\cap \mathring T^{pn}D'=\varnothing, \ n>0.$$
In fact, the above series diverges. Therefore, we apply our idea to tensor powers of Poisson suspensions.
First, let us verify that
$$\mu_\circ (\mathring RD\cap \mathring T^{-n}\mathring R \mathring T^n D)< const \cdot \mu(T^n U\cap U).$$
Let the ideal case occur when $\mu(T^n U\cap U)=0$ for $n>0$. Then $T^{-n}R T^n$ acts identically on $U$ and $\mathring T^{-n}\mathring R\mathring T^n D=D$. This implies the repulsion effect:
$$\mu_\circ ((\mathring R^{-1}\mathring T \mathring R)^n D\cap \mathring T^n D)=0, \ n>0.$$ If the intersection measures $T^n U\cap U$ decrease in a power-law manner, the repulsion will be provided by the tensor powers of the Poisson suspensions $\mathring S,
\mathring T$. Let us verify this. It follows from the definition of the Poisson measure that $$\mu_\circ (\mathring RD\cap \mathring T^{-n}\mathring R\mathring T^n D) < const \cdot \mu(T^n U\cap U).$$
The correlations $\mu(T^n U\cap U)$ for Sidon $T\in C(1)$ satisfy the condition $$\sum_n \mu(T^n U\cap U)^4\ <\ \infty,$$
which follows from the remark at the end of section \S 3 for $U \subset X_1$.
However, the Sidon automorphism of class $T\in C(0)$ is also suitable for this purpose. From the above we now obtain
$$\sum_n \mu_\circ^{\otimes 4}
((\mathring S^nD)^{\times 4}\cap (\mathring T^n D)^{\times 4})\ <\
(const)^4\  \ \sum_n \mu(T^n F\cap F)^4\ <\ \infty.$$
For a sufficiently large $p$ there exists $D'\subset D^{\times 4}$ of positive $\mu_\circ^{\otimes 4}$-measure such that
$$(\mathring S^{\times 4})^pD'\ \cap \ (\mathring T^{\times 4})^pD'=\varnothing.$$
Poisson suspensions over automorphisms of rank one and their tensor powers are known to have zero entropy \cite{Ja}. Theorem 3 is proved.

The automorphisms $\mathring S, \mathring T$ have mixed spectrum consisting of  singular component and  countably multiple Lebesgue component. The proposed method does not allow us to find examples of repulsive pairs $S,T$ with singular spectrum.

\vspace{3mm}

\section{Qestions}

\ \ \  1. It would be interesting to know which invariants of conjugate mixing automorphisms $S,T$ are incompatible with \it repulsibility: \rm  $S^nD\cap T^nD=\varnothing$ for all $n>0$ and some set $D$ of positive measure.

\vspace{5mm}
What we can say about the repulsibility  $S^{p(n)}D\cap T^{q(n)}D=\varnothing$?

\vspace{5mm}
2. Can some tensor products of infinite automorphisms have  simple Lebesgue \\ spectrum or at least  absolutely continuous spectra  of   finite multiplicity?\rm

\vspace{5mm}
3.  A joining   $\nu$ of actions $T_1,\dots,T_n$ is a probability $T_1\times \dots\times T_n$- invariant measure on $X^{\times n}$  with marginal  projections equal   $\mu$. If $T_1,\dots,T_n=T$  such a joining is called self-joining of $T$ of order $n$.
  We say that the action $T$ belongs to the class $PID$,  if every self-joining of order $n>2$ such that all projections onto the $2$-dimensional marginals are  $\mu^2$, is  $\mu^{\times n}$, the product of $n$ copies of $\mu$. Positive entropy and the presence of a discrete spectrum both are incompatible with $PID$, so 
Rudolph and   del Junco asked:  
  \do automorphisms of zero entropy and continuous spectra possess    $PID$? \rm
No answer for 40 years. There is Host's theorem:  singular spectrum  implies $PID$.  

\vspace{5mm}
Do ergodic Poisson suspensions  and Gaussian automorphisms with zero entropy  and Lebesgue component in spectrum  possess    $PID$? \rm

\vspace{25mm}
The author thanks Nikos Frantzikinakis and Jean-Paul Thouvenot for stimulating questions and useful discussions.

%$\end{fulltext}

\end{document}